\documentclass[a4paper,12pt,twoside]{article}
\usepackage[latin2]{inputenc}
\usepackage{amssymb}
\usepackage{amsmath}
\usepackage[mathscr]{eucal}
\usepackage[dvips]{graphicx}
\makeatletter

\def\Lifts{{\cal L}}
\def\Forb{\mathop{\mathrm{Forb_h}}\nolimits}

\def\UP{\mathop{\mathrm{UP}}\nolimits}
\def\Age{\mathop{\mathrm{Age}}\nolimits}
\def\Rel{\mathop{\mathrm{Rel}}\nolimits}
\def\CSP{\mathop{\mathrm{CSP}}\nolimits}

\def\Block{\mathop{\mathrm{Block}}\nolimits}
\def\ig{\mathop{\mathrm{IG}}\nolimits}
\def\sh{\mathop{\mathrm{Sh}}\nolimits}

\def\relsys#1{\mathbf {#1}}
\def\extsys#1{\mathbf #1'}
\def\rel#1#2{R_{\mathbf{#1}}^{#2}}
\def\ext#1#2{X_{\mathbf{#1}}^{#2}}
\def\extl#1#2{X_{#1}^{#2}}
\def\F{{\cal F}}
\def\K{{\cal K}}
\def\Piece{{\cal P}}
\def\APiece{{\cal A}}
\def\Incompatible{{\cal I}}

\def\Fraisse{Fra\"{\i}ss\' e}

\newtheorem{defn}{Definition}[section]

\newtheorem{corollary}{Corollary}[section]
\newtheorem{prop}{Proposition}[section]
\newtheorem{observation}{Observation}[section]

\newtheorem{thm}{Theorem}[section]
\newtheorem{lem}[thm]{Lemma} 

\DeclareRobustCommand{\qed}{%
  \ifmmode \mathqed
  \else
    \leavevmode\unskip\penalty9999 \hbox{}\nobreak\hfill
    \quad\hbox{$\square$}%
  \fi
}
\let\QED@stack\@empty
\let\qed@elt\relax
\newcommand{\pushQED}[1]{%
  \toks@{\qed@elt{#1}}\@temptokena\expandafter{\QED@stack}%
  \xdef\QED@stack{\the\toks@\the\@temptokena}%
}
\newcommand{\popQED}{%
  \begingroup\let\qed@elt\popQED@elt \QED@stack\relax\relax\endgroup
}
\def\popQED@elt#1#2\relax{#1\gdef\QED@stack{#2}}
\newcommand{\qedhere}{%
  \begingroup \let\mathqed\math@qedhere
    \let\qed@elt\setQED@elt \QED@stack\relax\relax \endgroup
}

\providecommand{\proofname}{Proof}
\newenvironment{proof}[1][\proofname]{\par
  \pushQED{\qed}%
  \normalfont \topsep6\p@\@plus6\p@\relax
  \trivlist
  \item[\hskip\labelsep \bf {#1\ignorespaces.}]\ignorespaces
}{%
\popQED\endtrivlist
\par
}

\begin{document}

\title{Universal structures with forbidden homomorphisms}

\author{
   \normalsize {\bf Jan Hubi\v cka$^*$}\\
   \normalsize {\bf Jaroslav Ne\v set\v ril}\thanks {The Computer Science Institute of Charles University (IUUK) is supported by grant ERC-CZ LL-1201 of the Czech Ministry of Education and CE-ITI P202/16/6061 of GA\v CR.}
\\
{\small Computer Science Institute of Charles University (IUUK)}\\
   {\small Charles University}\\
   {\small Malostransk\' e n\' am. 25, 11800 Praha, Czech Republic}\\
   {\normalsize 118 00 Praha 1}\\
   {\normalsize Czech Republic}\\
   {\small \{hubicka,nesetril\}@iuuk.mff.cuni.cz}
}

\date{}
\maketitle
\begin{abstract}
We relate the existence problem of universal objects to the properties of
corresponding enriched categories (lifts or expansions). In particular,
extending earlier results, we prove that for every (possibly infinite) regular set $\F$ of finite
connected structures there exists a (countable) $\omega$-categorical universal
structure $\relsys{U}$ for the class $\Forb(\F)$ (of all countable structures
not containing any homomorphic image of a member of $\F$).  We employ a
technique known as homogenization:  The universal object $\relsys U$ is the
shadow (reduct) of an ultrahomogeneous structure $\relsys U'$. 

We also put the results of this paper in the context of homomorphism dualities
and constraint satisfaction problems. This leads to an alternative proof of the
characterization of finite dualities (given by Tardif and Ne\v set\v ril) as
well as of the characterization of infinite-finite dualities for classes of
relational trees given by P.~L.~Erd\H{o}s, P\'alv\" olgyi, Tardif and Tardos.

The notion of regular families of structures is motivated by the recent
characterization of infinite-finite dualities for classes of relational forests
(itself related to regular languages). We show how the notion of a regular
family of relational trees can be extended to regular families of relational
structures.  This gives a partial characterization of the existence of a
(countable) $\omega$-categorical universal object for classes $\Forb(\F)$.
\end{abstract}
\section {Introduction}
We review first a few well known concepts and facts.

A {\em relational structure} (or simply {\em structure}) $\relsys{A}$ is a pair
$(A,(\rel{A}{i}:i\in I))$, where $\rel{A}{i}\subseteq A^{\delta_i}$ (i.e.,
$\rel{A}{i}$ is a $\delta_i$-ary relation on $A$). The family $(\delta_i: i\in
I)$ is called the {\em type} $\Delta$. The type is usually fixed and understood
from the context.  We consider only finite types.  If the set $A$ is finite we
call $\relsys A$ a {\em finite structure}.  We consider only countable or
finite structures.  The class of all (countable) relational structures of type
$\Delta$ will be denoted by $\Rel(\Delta)$.  The class $\Rel(\Delta),
\Delta=(\delta_i; i\in I)$, is fixed throughout this paper. Unless otherwise
stated all structures $\relsys{A}, \relsys{B},\ldots$ belong to $\Rel(\Delta)$. 

A {\em homomorphism} $f:\relsys{A}\to \relsys{B}=(B,(\rel{B}{i}:i\in I))$ is a
mapping $f:A\to B$ such that $(x_1,x_2,\ldots, x_{\delta_i})\in \rel{A}{i}$
implies $(f(x_1),f(x_2),\ldots,f(x_{\delta_i}))\in \rel{B}{i}$, for each $i\in
I$.  For given structures $\relsys{A}$ and $\relsys{B}$ we will denote the
existence of homomorphism $f:\relsys{A}\to \relsys{B}$ by
$\relsys{A}\to\relsys{B}$ and the non-existence by $\relsys{A}\nrightarrow
\relsys{B}$.  If $f$ is one-to-one then $f$ is called a {\em monomorphism}. A
monomorphism $f$ such $(x_1,x_2,\ldots, x_{\delta_i})\in \rel{A}{i}$ if and
only if $(f(x_1),f(x_2),\ldots,f(x_{\delta_i}))\in \rel{B}{i}$ for each $i\in
I$ is called an {\em embedding}. 

Given a family of relational structures $\F$, by $\Forb(\F)$ we denote the
class of all relational structures $\relsys{A}$ for which there is no
homomorphism $\relsys{F}\to \relsys{A}$, for any $\relsys{F}\in \F$. Formally,
$$\Forb(\F)= \{\relsys{A}; \forall_{\relsys{F}\in \F}\relsys{F}\nrightarrow
\relsys{A}\}.$$

Given a class $\K$ of countable structures, an object $\relsys{U}\in\K$ is
called {\em universal} for $\K$ if for every object
$\relsys{A}\in \K$ there exists an embedding $\relsys{A}\to \relsys{U}$.

For a class $\K$ of countable relational structures, we denote by $\Age(\K)$
the class of all finite structures isomorphic to a substructure of some
$\relsys{A}\in \K$ and call it the {\em age of $\K$}. Similarly, for a
relational structure $\relsys{A}$, the age of $\relsys{A}$, $\Age(\relsys{A})$,
is $\Age(\{\relsys{A}\})$.

A structure $\relsys A$ is {\em ultrahomogeneous} (sometimes also simply called
{\em homogeneous}) if every isomorphism between two induced finite
substructures of $\relsys A$ can be extended to an automorphism of $\relsys A$. 
A structure $\relsys{G}$ is {\em generic} for the class $\K$ if it is universal
for $\K$ and ultrahomogeneous.

The key property of the age of any ultrahomogeneous structure is described by
the following concept.
\begin{figure}
\centerline{\includegraphics{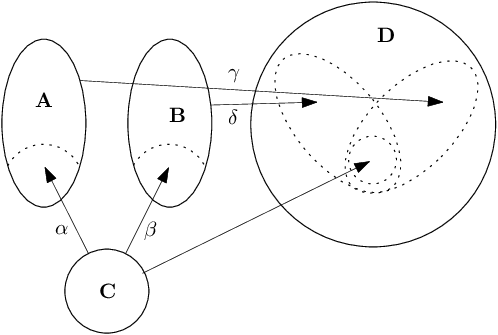}}
\caption{Amalgamation of $(\relsys{A},\relsys{B},\relsys{C}, \alpha, \beta)$.}
\label{amalgamfig0}
\end{figure}
\label{amalgamationclassdef}
Let $\relsys{A},\relsys{B},\relsys{C}$ be relational structures, $\alpha$ an
embedding of $\relsys{C}$ into $\relsys{A}$, and $\beta$ an embedding of
$\relsys{C}$ into $\relsys{B}$.  An {\em amalgamation of $(\relsys{A},
\relsys{B}, \relsys{C}, \alpha, \beta)$} is a triple
$(\relsys{D},\gamma,\delta)$, where $\relsys{D}$ is a relational structure,
$\gamma$ an embedding $\relsys{A}\to \relsys{D}$ and $\delta$ an embedding
$\relsys{B}\to\relsys{D}$ such that $\gamma\circ\alpha = \delta\circ\beta$.
Less formally, an amalgamation ``glues together'' the structures $\relsys{A}$
and $\relsys{B}$ into a single substructure of $\relsys{D}$ such that copies of
$\relsys{C}$ coincide.  See Figure~\ref{amalgamfig0}.

Often the vertex sets of structures $\relsys{A}$, $\relsys{B}$ and $\relsys{C}$
can be chosen in such a way that the embeddings $\alpha$ and $\beta$ are
identity mappings.  In this case, for brevity, we shall call an amalgamation of
$(\relsys{A},\relsys{B},\relsys{C},\alpha,\beta)$ simply an {\em amalgamation
of $\relsys{A}$ and $\relsys{B}$ over $\relsys{C}$}.  Similarly, for an
amalgamation $(\relsys{D}, \gamma, \delta)$ of a given $(\relsys{A},
\relsys{B}, \relsys{C}, \alpha, \beta)$ we are often interested in the
structure $\relsys{D}$ alone.  In this case we shall call the structure
$\relsys{D}$ an amalgamation of $(\relsys{A}, \relsys{B}, \relsys{C}, \alpha,
\beta)$ (omitting the embeddings $\gamma$ and $\delta$).

We say that an amalgamation is {\em strong} when $$\gamma(A)\cap \delta (B)=\gamma (\alpha (C)) = \delta (\beta (C)).$$
Less formally, a strong
amalgamation glues together $\relsys{A}$ and $\relsys{B}$ with an overlap no
greater than the copy of $\relsys{C}$ itself.  A strong amalgamation is {\em
free} if there are no relations of $\relsys{D}$ spanning both vertices of
$\gamma(A)$ and $\delta(B)$ that are not images of some relations of structure
$\relsys{A}$ or $\relsys{B}$ via the embedding $\gamma$ or $\delta$,
respectively.

A class $\K$ of finite relational structures is called an {\em amalgamation
class} if the following conditions hold:
\begin{enumerate}
\item ({\em Hereditary property}) For every $\relsys{A}\in \K$ and induced substructure $\relsys{B}$ of $\relsys{A}$ we have $\relsys{B}\in \K$.
\item ({\em Amalgamation property}) 
For $\relsys{A},\relsys{B},\relsys{C}\in \K$ and $\alpha$ an embedding of
$\relsys{C}$ into $\relsys{A}$, $\beta$ an embedding of $\relsys{C}$ into
$\relsys{B}$, there exists $(\relsys{D},\gamma,\delta), \relsys{D}\in \K$, that is an amalgamation of
$(\relsys{A}, \relsys{B}, \relsys{C}, \alpha, \beta)$.
\item $\K$ is closed under isomorphism.
\item $\K$ has only countably many mutually non-isomorphic structures. 
(This is always the case in our setting of finite types).
\end{enumerate}

The following classical result establishes the correspondence between
amalgamation classes and ultrahomogeneous structures.
\begin{thm}[\Fraisse{} \cite{F, Hoges}]
\label{fraissethm}
(a) A class $\K$ of finite structures is the age of a countable
ultrahomogeneous structure $\relsys{G}$ if and only if $\K$ is an amalgamation
class. 

(b) If the conditions of (a) are satisfied then the structure $\relsys{G}$ is
unique up to isomorphism. 
\end{thm}

The ultrahomogeneous structure $\relsys{G}$ such that $\Age(\relsys{G})=\K$ is
called the {\em \Fraisse{} limit} of $\K$.  We say that structure $\relsys{A}$
is {\em younger} than structure $\relsys{B}$ if $\Age(\relsys{A})$ is a subset
of $\Age(\relsys{B})$.  Every ultrahomogeneous structure $\relsys{G}$ has the
property that it is universal (and also generic) for the class $\K$ of all countable structures
younger than $\relsys{G}$.

A countably infinite structure is called {\em $\omega$-categorical} if all
countable models of its first order theory are isomorphic.  We use the
following characterization of $\omega$-categorical structures given by Engeler
\cite{Engeler}, Ryll-Nardzewski \cite{Nardzewsky} and Svenonius
\cite{Svenonius}.

\begin{thm}
\label{Nardzewsky}
For a countable first order structure $\relsys{A}$, the following conditions
are equivalent:
\begin{enumerate}
\item $\relsys{A}$ is $\omega$-categorical.
\item The automorphism group of $\relsys{A}$ has only finitely many orbits on
$n$-tuples, for every $n$.
\end{enumerate}
\end{thm}

\paragraph{Lifts and shadows.}

Let $\Delta'=(\delta'_i;i\in I')$ be a type containing type $\Delta$. (By this
we mean $I\subseteq I'$ and $\delta'_i=\delta_i$ for $i\in I$.) Then every
structure $\relsys{X}\in \Rel(\Delta')$ may be viewed as a structure
$\relsys{A}=(A,(\rel{A}{i}; i\in I))\in \Rel(\Delta)$ together with some
additional relations for $i\in I'\setminus I$. To make this more
explicit, these additional relations will be denoted by $\ext{X}{i}, i\in
I'\setminus I$. Thus a structure $\relsys{X}\in \Rel(\Delta')$ will be written
as $$\relsys{X}=(A,(\rel{A}{i};i\in I),(\ext{X}{i};i\in I'\setminus I)),$$ and,
abusing notation, more briefly as
$$\relsys{X}=(\relsys{A},\ext{X}{1},\ext{X}{2},\ldots, \ext{X}{N}).$$

We call $\relsys{X}$ a {\em lift} of $\relsys{A}$ and $\relsys{A}$ is called
the {\em shadow} of $\relsys{X}$. In this sense the class $\Rel(\Delta')$ is
the class of all lifts of $\Rel(\Delta)$.  Conversely, $\Rel(\Delta)$ is the
class of all shadows of $\Rel(\Delta')$. If all extended relations are unary,
the lift is called {\em monadic}.  In the context of monadic lifts, the {\em
color} of vertex $v$ is the set $\{i;(v)\in \ext{U}{i}\}$.
Note that a lift is also in the model-theoretic setting called an {\em
expansion} (as we are expanding our relational language) and a shadow a {\em reduct} (as we are reducing it).  (Our terminology is motivated by a
computer science context; see \cite{KunN}.) Unless stated explicitely, we shall
use letters $\relsys{A}, \relsys{B}, \relsys{C}, \ldots$ for shadows (in
$\Rel(\Delta)$) and letters $\relsys{X}, \relsys{Y}, \relsys{Z}$ for lifts (in
$\Rel(\Delta'))$.

For a lift $\relsys X=(\relsys{A}, \ext{X}{1},\ldots,\ext{X}{N})$ we denote by
$\sh(\relsys{X})$ the relational structure $\relsys{A}$, i.e. its shadow.
($\sh$ is called a {\em forgetful functor}.) Similarly, for a class $\K'$ of
lifted objects, we denote by $\sh(\K')$ the class of all shadows of structures
in $\K'$.

\paragraph{Homogenization.}
Many naturally defined classes $\K$ of relational structures contain universal
structures that are $\omega$-categorical.  Because $\omega$-ca\-te\-go\-ri\-city can be
seen as a weaker notion of ultrahomogeneity, it is natural to construct
$\omega$-categorical universal structures as shadows of ultrahomogeneous
structures. Such a construction is called {\em homogenization}.  Covington
\cite{Covington} provided a sufficient condition for the existence of a
universal structure for a given class $\K$ that is a shadow of an
ultrahomogeneous structure by means of amalgamation failures. This concept in
fact relaxes the \Fraisse{} Theorem.

However, not all universal structures are constructed by means of
homogenization.  A necessary and sufficient condition for the existence of a
universal structure for the class defined by forbidden monomorphisms from a
finite family $\F$ of connected graphs was given by Cherlin, Shelah and Shi
\cite{CherlinShelahShi}.  Here the classes are characterized by means of local
finiteness of the algebraic closure operator.  The techniques of \cite{CherlinShelahShi} 
are motivated by proofs of the non-existence of a universal structure for
a given class. The universal structure is not constructed by an explicit
amalgamation argument.

\paragraph{Our motivation and results.}
Our motivation stems from several sources.  First, we seek a more streamlined
and combinatorial proof of the following corollary of the aforementioned result
of Cherlin, Shelah and Shi.

\begin{thm}[\cite{CherlinShelahShi}]
\label{homclosed}
For every finite family $\F$ of finite connected graphs there is an
$\omega$-categorical universal graph for the class $\Forb(\F)$.
\end{thm}

We prove a stronger form of Theorem~\ref{homclosed} by an explicit amalgamation
argument. A similar construction can be also be obtained by characterizing the
amalgamation failures and applying homogenization method.  Our
lifts are, however, different and more effective.

We are interested in the structure of lifts constructed for a given (possibly
infinite) family $\F$.  In special cases we relate lifts to the concept of
homomorphism dualities. Motivated by a recent characterization of
infinite-finite dualities by P. L. Erd\H{o}s, P\'alv\" olgyi, Tardif, Tardos
\cite{Tardif}, we introduce a notion of regular families of relational
structures. These (possibly infinite) families of structures generalize regular
forests, used in \cite{Tardif} to characterize infinite-finite dualities.  In
fact the regular families of trees (and forests) corresponds to well
established notion of recognizable tree languages, see \cite{tata}.

In Section~\ref{sec:mainthm} we strengthen Theorem~\ref{homclosed} by proving
the existence of a universal structure for $\Forb(\F)$, where $\F$ is a regular
family of finite connected structures.

In Section~\ref{sec:classification} we show the non-existence of an
$\omega$-categorical universal structure for $\Forb(\F)$ for certain
non-regular families $\F$, and give a partial characterization of such families.

Finally, in Section~\ref{sec:dualities} we relate our results to homomorphism
dualities and constraint satisfaction problems.  We show that for the classes
$\F$ consisting of regular relational trees the universal structure has a
finite retract. This gives an alternative construction of graph duals and also
an alternative proof of the characterization of homomorphism dualities.

\section{Regular families of structures and $\F$-lifts}
\label{sec:regular}
Let $\F$ be a fixed family of finite connected relational structures.  For the
construction of a universal structure of $\Forb(\F)$ we use special lifts,
called $\F$-lifts.  The definition of an $\F$-lift is easy and resembles
decomposition techniques standard in graph theory, and thus we adopt a similar
terminology.  First we review some elementary graph-theoretic notions, see
\cite{MatN, Bolobas} for details.

For a structure $\relsys{A}=(A,(\rel{A}{i},i\in I))$, the {\em Gaifman graph}
(in combinatorics often called {\em 2-section}) is the graph $G_\relsys{A}$
with vertices $A$ and all those edges which are a subset of a tuple of a
relation of $\relsys{A}$, i.e., $G=(A,E),$ where ${x,y}\in E$ if and only if
$x\neq y$ and there exists a tuple $\vec{v}\in \rel{A}{i}$, $i\in I$, such that
$x,y\in \vec{v}$.

We adopt the following standard graph-theoretic notions for relational structures.
We call a structure $\relsys{A}$ {\em connected} if its Gaifman graph
$G_\relsys{A}$ is connected.
For a structure $\relsys{A}$ and subset of its vertices $B\subseteq A$, we
denote by $N_\relsys{A}(B)$ the {\em neighborhood} of the set $B$, that is all
vertices of $A\setminus B$ connected in the Gaifman graph $G_\relsys{A}$ by an edge
to a vertex of $B$.
We denote by $\relsys{A}\setminus B$ the structure induced on $A\setminus B$ by $\relsys{A}$.
 Similarly we denote by $G_\relsys{A}\setminus B$ the graph
created from the Gaifman graph $G_\relsys{A}$ by removing the vertices in $B$.

A {\em g-cut} in $\relsys{A}$ is a subset $C$ of $A$ such that the Gaifman
graph $G_\relsys{A}$ is disconnected by removing the set $C$. That is, there are
vertices $u,v\in A\setminus C$ that belong to the same connected component of
$G_\relsys{A}$ but to different connected components of
$G_\relsys{A}\setminus C$. A {\em cut} in $\relsys{A}$ is subset $C$ of $A$
such that there are vertices $u,v\in A\setminus C$ that belong to the same
connected component of $\relsys{A}$ but to different connected components of $\relsys{A}\setminus C$.

Observe that not every cut is a g-cut. With relations of arity greater than 2,
$G_{\relsys{A}\setminus C}$ may be different from $G_{\relsys{A}}\setminus C$.

For a g-cut $C$ in the relational structure $\relsys{A}$, a structure
$\relsys{A}_1$ is a {\em g-component} of $\relsys{A}$ with g-cut $C$ if
$\relsys{A}_1$ is induced by $\relsys{A}$ on some connected component of
$G_\relsys{A}\setminus C$.

We will make use of the following simple observation about the neighborhood and
g-components.
\begin{observation}
\label{cuty}
Let $\relsys{A}_1$ be a g-component of $\relsys{A}$ with g-cut $C$. Then the neighborhood $N_\relsys{A}(A_1)$ is a subset of $C$.
Moreover $N_\relsys{A}(A_1)$ is a g-cut and $\relsys{A}_1$ is one of the g-components of $\relsys{A}$ with g-cut $N_\relsys{A}(A_1)$.
\end{observation}

Given a structure $\relsys{A}$ with g-cut $C$ and two (induced) substructures
$\relsys{A}_1$ and $\relsys{A}_2$, we say that $C$ {\em g-separates}
$\relsys{A}_1$ and $\relsys{A}_2$ if there are g-components $\relsys{A}'_1\neq
\relsys{A}'_2$ of $\relsys{A}$ with g-cut $C$ such that $A_1\subseteq A'_1$ 
and $A_2\subseteq A'_2$.

\begin{defn}
Let $C$ be a g-cut in structure $\relsys{A}$. Let $\relsys{A}_1\neq \relsys{A}_2$
be two g-components of $\relsys{A}$ with g-cut $C$.  We call $C$ {\em
minimal g-separating g-cut} for $\relsys{A}_1$ and $\relsys{A}_2$ in $\relsys{A}$ if
$C=N_\relsys{A}(A_1)=N_\relsys{A}(A_2)$.
\end{defn}
For brevity, we can omit one or both $g$-components when speaking about
a minimal g-separating g-cut. Explicitly,
we call a g-cut $C$ {\em minimal g-separating} for $\relsys{A}_1$
in $\relsys{A}$ if there exists another structure $\relsys{B}$ such that $C$
is minimal g-separating for $\relsys{A}_1$ and $\relsys{B}$ in $\relsys{A}$.  A g-cut
$C$ is {\em minimal g-separating} in $\relsys{A}$ if there exists structures
$\relsys{B}_1$ and $\relsys{B}_2$ such that $C$ is minimal g-separating for
$\relsys{B}_1$ and $\relsys{B}_2$ in $\relsys{A}$.

The name of minimal g-separating g-cut is justified by the following (probably folkloristic) proposition.
\begin{prop}
Let $\relsys{A}$ be a connected relational structure, $C$ a g-cut in $\relsys{A}$ and
$\relsys{A}_1$ and $\relsys{A}_2$ (induced) substructures of $\relsys{A}$
g-separated by $C$. Then there exists a minimal g-separating g-cut $C'\subseteq C$ that
g-separates $\relsys{A}_1$ and $\relsys{A}_2$ in $\relsys{A}$.
Moreover if $N_\relsys{A}(A_1)\subseteq C$ (or, equivalently, $\relsys{A}_1$ is a g-component of $\relsys{A}$ with g-cut $C$), then $C'\subseteq N_\relsys{A}(A_1)$.
\end{prop}
\begin{proof}
\begin{figure}
\centerline{\includegraphics{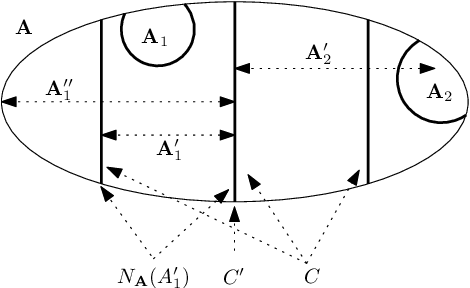}}
\caption{Construction of a minimal g-separating g-cut $C'$ separating $\relsys{A}_1$ and $\relsys{A}_2$ in $\relsys{A}$.}
\label{constructioncut}
\end{figure}
\label{prop:sep}
We will construct a series of g-cuts and g-components as depicted in Figure~\ref{constructioncut}.

Denote by $\relsys{A}'_1$ the g-component of $\relsys{A}$ with g-cut $C$
containing $\relsys{A}_1$ (and thus not containing $\relsys{A}_2$).  By Observation~\ref{cuty},
$N_\relsys{A}(A'_1)\subseteq C$ is a g-cut that g-separates $\relsys{A}'_1$ and
$\relsys{A}_2$ (because $\relsys{A}'_1$ is also g-component of $\relsys{A}$ with cut
$N_\relsys{A}(A'_1)$ and $\relsys{A}'_1$ does not contain $\relsys{A}_2$).  

Now
consider g-component $\relsys{A}'_2$ of $\relsys{A}$ with g-cut
$N_\relsys{A}(A'_1)$ containing $\relsys{A}_2$.  Put $C'=N_\relsys{A}(A'_2)$. By
Observation~\ref{cuty}, $C'\subseteq N_\relsys{A}(A'_1)\subseteq C$ is g-cut and
$\relsys{A}'_2$ (not containing $\relsys{A}_1$) is one of its g-components.

Denote by $\relsys{A}''_1$ the g-component of $\relsys{A}$ with cut $C'$
containing $\relsys{A}_1$. It follows that $C'$ g-separates $\relsys{A}''_1$ (that contains $\relsys{A}_1$) and $\relsys{A}'_2$ (that contains $\relsys{A}_2$).

To see that $C'$ is minimal g-separating for $\relsys{A}''_1$ 
and $\relsys{A}'_2$ it remains
to show that every vertex in $C'=N_\relsys{A}(A'_2)$ is also in $N_\relsys{A}(A''_1)$.  This is true
because every vertex of $C'$ is in $N_\relsys{A}(A'_1)$ and $\relsys{A}'_1$ is
substructure of $\relsys{A}''_1$. 

\end{proof}
Observe that every inclusion minimal g-cut is also minimal g-separating, but
not vice versa. Every minimal g-separating g-cut $C'\subset C$ that g-separates
$\relsys{A}_1$ and $\relsys{A}_2$ is however also inclusion minimal g-cut that
separates $\relsys{A}_1$ and $\relsys{A}_2$.  

If $C$ is a set of vertices then $\overrightarrow{C}$ will denote a tuple (of
length $|C|$) of all the elements of $C$. Alternatively, $\overrightarrow{C}$
is an arbitrary linear ordering of $C$.
A {\em rooted structure} $\Piece$ is a pair $(\relsys{P},\overrightarrow{R})$
where $\relsys{P}$ is a relational structure and $\overrightarrow{R}$ is a
tuple consisting of distinct vertices of $\relsys{P}$. $\overrightarrow{R}$ is
called the {\em root} of $\Piece$ and the size of $\overrightarrow{R}$ is the
{\em width} of $\Piece$.  We say that rooted structures
$\Piece_1=(\relsys{P}_1,\overrightarrow{R}_1)$ and
$\Piece_2=(\relsys{P}_2,\overrightarrow{R}_2)$ are {\em isomorphic} if there is
a function $f:P_1\to P_2$ that is an isomorphism of structures $\relsys{P}_1$
and $\relsys{P}_2$ and $f$ restricted to $\overrightarrow{R}_1$ is a monotone
bijection between $\overrightarrow{R}_1$ and $\overrightarrow{R}_2$ (we denote
this $f(\overrightarrow{R}_1)=\overrightarrow{R}_2$).

The following is the principal notion of this paper:
\begin{defn}
Let $\relsys{A}$ be a connected relational structure and $R$ a minimal g-separating
g-cut for g-component $\relsys{C}$ in $\relsys{A}$. A {\em piece} of a relational structure
$\relsys A$ is then a rooted structure $\Piece=(\relsys{P},\overrightarrow{R})$, where the
tuple $\overrightarrow{R}$ consists of the vertices of the g-cut $R$ in a (fixed)
linear order and $\relsys {R}$ is a structure induced by $\relsys{A}$ on $C\cup R$.
\end{defn}
Note that since $\relsys{P}$ is the union of a g-component and its neighborhood it follows that
the pieces of a connected structure are always connected structures.
\begin{figure}
\centerline{\includegraphics{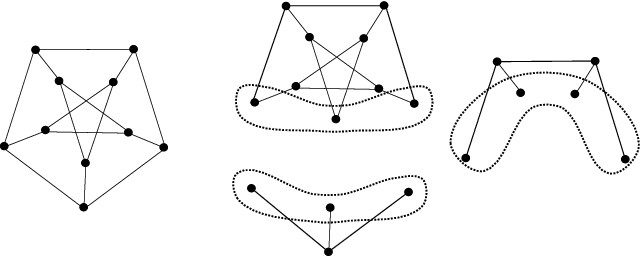}}
\caption{Pieces of the Petersen graph up to isomorphism (and permutations of roots).
G-cut of size 4 splits the graph into 3 isomorphic pieces.}
\label{Petersoni}
\end{figure}
As an example, pieces of the Petersen graph are shown in Figure~\ref{Petersoni}.

Given rooted structures $(\relsys{P},\overrightarrow{R})$ and
$(\relsys{P}',\overrightarrow{R}')$ such that $|R|=|R'|$, denote by
$(\relsys{P},\overrightarrow{R})\oplus (\relsys{P}',\overrightarrow{R}')$ the
(possibly rooted) structure created as a free amalgam of $\relsys{P}$ and $\relsys{P}'$
with corresponding roots being identified (in the order of $\overrightarrow{R}$ and $\overrightarrow{R}'$).
Note that $(\relsys{P},\overrightarrow{R})\oplus (\relsys{P}',\overrightarrow{R}')$ is defined only
if the rooted structure induced by $\relsys{P}$ on $\overrightarrow{R}$ is isomorphic to the rooted
structure induced by $\relsys{P}'$ on $\overrightarrow{R}'$.

Recall that $\F$ is a fixed family of finite connected relational structures.
A piece $\Piece=(\relsys{P},\overrightarrow{R})$ is {\em incompatible} with
a rooted structure $\APiece$ if $\Piece \oplus \APiece$ is defined and there exists $\relsys{F}\in \F$
that is isomorphic to $\Piece \oplus \APiece$.  
(In other words, there exists $\relsys{F}'$ isomorphic to some $\relsys{F}''\in \F$, such that $\Piece$ is a piece of $\relsys{F}'$ and $\APiece$ is the structure induced on $F'\setminus (P\setminus R)$
by $\relsys{F}'$ rooted by $\overrightarrow{R}$.)

Assign to each piece $\Piece$ a set $\Incompatible_\Piece$ containing all rooted
structures that are incompatible with $\Piece$.  For two pieces $\Piece_1$ and
$\Piece_2$ put $\Piece_1 \sim \Piece_2$ if and only if
$\Incompatible_{\Piece_1}=\Incompatible_{\Piece_2}$.  Observe that every
equivalence class of $\sim$ contains pieces of the same width $n$. We also
call $n$ the {\em width} of the equivalence class of $\sim$.

\begin{defn}
A family of finite structures $\F$ is called {\em regular} if 
there are only finitely many equivalence classes of $\sim$ on the family of all
pieces of structures in family $\F$.
\end{defn}
The notion of regular family is a generalization of that of a regular family of
forests, introduced in \cite{Tardif}. The term used in \cite{Tardif} was motivated by the connection
to regular languages. This is explained in the following examples.

\paragraph{Example.}
All finite families $\F$ of finite structures are regular. Examples of infinite
families $\F$ include the following:
\begin{enumerate}
\item The family $\F_{\mathrm {odd}}$ consisting of all graph cycles of odd
length. All pieces of $\F_{\mathrm {odd}}$ are paths rooted by initial vertex and terminal vertex.
There are only two equivalence classes of the pieces: paths of odd length and
paths of even length.
\item The family $\F_{\mathrm {oriented}}$ consisting of those orientations of
graph cycles where all edges are oriented in the same direction.  Pieces of
$\F_{\mathrm {oriented}}$ are oriented paths with all edges in a forward direction with roots on
initial and terminal vertex.  Consequently there are only two equivalence
classes of pieces: paths with first root on initial vertex and second root on
terminal vertex, and paths with first root on terminal vertex and second root
on initial vertex.

\item Oriented paths can be described by words on alphabet
$\{\leftarrow,\rightarrow\}$.  It follows that every language of words on this
alphabet corresponds to a family of oriented paths.  It is not difficult to
show that all regular languages correspond to a regular family of paths.
Consequently regular families may have a rich structure; see \cite{Tardif}.

Consider for example the family created by words of the form
$\rightarrow\rightarrow(\rightarrow\leftarrow\rightarrow)^n\rightarrow\rightarrow$,
$n\geq 1$, where $(\rightarrow\leftarrow\rightarrow)^n$ stands for $n$
repetitions of $\rightarrow\leftarrow\rightarrow$.
All these paths are cores and form an antichain. Several other examples of
regular families of directed graphs are discussed in
\cite{Tardif2}.
\end{enumerate}

We continue our construction with the following:
\begin{defn}
We denote by $E_1,\ldots E_N$ the equivalence classes of $\sim$ corresponding
to pieces of structures in $\F$. Put $I'=\{1,2,\ldots, N\}$.  A relational
structure $\relsys X=(\relsys{A}, (\ext{X}{i}:i\in I'))$ is called an {\em
$\F$-lift} of the relational structure $\relsys A$ when the arities of
the relations $\ext{X}{i}, i\in I'$, correspond to the width of $E_i$.

For a relational structure $\relsys A$, we define the {\em canonical lift}
$$L({\relsys A})=(\relsys{A},\extl{L(\relsys{A})}{1},\extl{L(\relsys{A})}{2},\ldots, \extl{L(\relsys{A})}{N})$$ by putting $(v_1,v_2,\ldots,v_l)\in \extl{L(\relsys{A})}{i}$ if
and only if there is a piece $\Piece=(\relsys{P},\overrightarrow{R})\in E_i$
such that there is a homomorphism $f:\relsys{P}\to\relsys{A}$
with $f(\overrightarrow {R})=(v_1,v_2,\ldots,v_l)$.
\end{defn}

\paragraph{Example.}
As an illustration we provide an explicit description of some canonical lifts of the
regular families discussed above.
\begin{enumerate}
\item For the family $\F_{\mathrm {odd}}$ there are two new binary relations.
In a canonical lift $L(\relsys{A})$ there is $(u,v)\in \extl{L(\relsys{A})}{1}$ if and
only if there is a walk of odd length between vertices $u$ and $v$ and
$(u,v)\in \extl{L(\relsys{A})}{2}$ if and only if there is a walk of even length between
$u$ and $v$.

For $\relsys{A}\in \Forb(\F_{\mathrm {odd}})$ there is no $(u,v)$ such that
$(u,v)\in \extl{L(\relsys{A})}{1}$ and $(u,v)\in \extl{L(\relsys{A})}{2}$.  This means that odd
cycles can be recognized by the existence of both a walk of even length and a
walk of odd length in between a given pair of vertices.

\item For the family $\F_{\mathrm {oriented}}$ there are two new binary relations. In a canonical lift $L(\relsys{A})$ there is $(u,v)\in \extl{L(\relsys{A})}{1}$ if and only if
there is an oriented walk from  $u$ to $v$ and $(u,v)\in \extl{L(\relsys{A})}{2}$ if and only if there is an oriented walk from $v$ to $u$.

For $\relsys{A}\in \Forb(\F_{\mathrm {oriented}})$ there is no $(u,v)$ such
that both $(u,v)\in \extl{L(\relsys{A})}{1}$ and $(v,u)\in \extl{L(\relsys{A})}{1}$.  The same
holds for $\extl{L(\relsys{A})}{2}$ and in fact the second relation is fully redundant in
our construction and can be ignored.
\end{enumerate}

\section{Construction of universal structures}
\label{sec:mainthm}
\begin{thm}
\label{mainthm}
Let $\F$ be a regular family of finite connected relational structures (of a
finite type).  Then there exists an ultrahomogeneous lift $\extsys{U}$ with only finitely
many new relations such that its shadow $\sh(\extsys{U})$ is a universal
structure for the class $\Forb(\F)$.

Moreover, the lift $\extsys{U}$ can be constructed in the following way.
Denote by $\Lifts$ the class of all induced substructures (sublifts) of
canonical lifts $L(\relsys{A})$, $\relsys{A}\in\Forb(\F)$. Then $\Age(\Lifts)$
is a amalgamation class (closed for strong amalgams whose shadows are free
amalgams) and $\extsys{U}$ is the \Fraisse{} limit of $\Age(\Lifts)$.

Denote by $n$ the maximal size of a minimal g-separating g-cut in a structure
in $\F$ (by regularity of $\F$ the size of g-cuts is bounded). Then the arity
of extended relations is bounded by $n$.
\end{thm}

This theorem will be proved in the rest of this section. We take time for a simple Lemma.

Given a piece $\Piece=(\relsys{P},\overrightarrow{R})$ of structure
$\relsys{F}$, we call $\Piece'=(\relsys{P}',\overrightarrow{R}')$ a {\em
subpiece} of $\Piece$ if $\Piece'$ is a piece of $\relsys{F}$, $P'\subset P$.
We show that a subpiece can be freely replaced by an
equivalent subpiece without changing the equivalence class of a given piece.
\begin{lem}
Let $\Piece_1=(\relsys{P}_1,\overrightarrow{R}_1)$ be a piece of structure
$\relsys F_1\in \F$, $\Piece_1'=(\relsys{P}_1',\overrightarrow{R}'_1)$ be a
subpiece of $\Piece_1$ and $\Piece'_2=(\relsys{P}'_2,\overrightarrow{R}'_2)$ a piece such that $\Piece'_1\sim \Piece'_2$.

Create $\Piece_2=(\relsys{P}_2,\overrightarrow{R}_2)$ as a copy of $\Piece_1$ with
$\relsys{P}'_1$ replaced by $\relsys{P}'_2$ identifying $\overrightarrow{R}'_1$ with $\overrightarrow{R}'_2$.  Then
$\Piece_2$ is an isomorphic copy of a piece of some $\relsys F_2\in \F$, and
moreover $\Piece_1\sim \Piece_2$.
\label{eqkousky}
\end{lem}
\begin{proof}
\begin{figure}
\centerline{\includegraphics{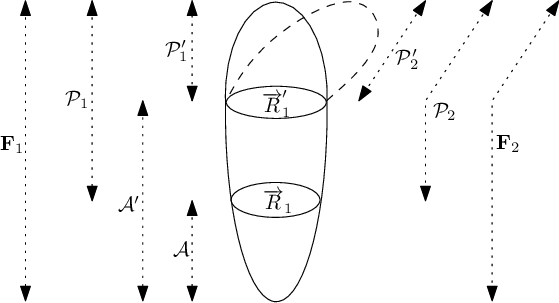}}
\caption{Replacing subpieces.}
\label{fig:replacing_pieces}
\end{figure}
See Figure \ref{fig:replacing_pieces}.
Consider some $\APiece\in\Incompatible_{\Piece_1}$. By definition
$\relsys{F}_1=\Piece_1\oplus\APiece$ is isomorphic to some structure $\F$.
$\Piece'_1$ is also a piece of $\relsys{F}_1$ and thus there exists a rooted structure $\APiece'$ such that $\Piece'_1\oplus\APiece'=\relsys{F}_1$.
It follows that $\APiece'\in \Incompatible_{\Piece'_1}$.
Because $\Piece'_2\sim \Piece'_1$ and thus $\Incompatible_{\Piece'_2}=\Incompatible_{\Piece'_1}$,
we also know that $\relsys{F}_2=\Piece'_2\oplus \APiece'$ is isomorphic to some structure in $\F$.
$\Piece_2$ is a copy of $\Piece_1$ with $\relsys{P}'_1$ replaced by $\relsys{P}'_2$.
By construction of $\APiece'$ we have $\Piece_2\oplus\APiece=\Piece'_2\oplus \APiece'=\relsys{F}_2$ and consequently $\APiece\in \Incompatible_{\Piece_2}$.
We get $\Incompatible_{\Piece_1}\subseteq \Incompatible_{\Piece_2}$.  By symmetry
we also have $\Incompatible_{\Piece_2}\subseteq \Incompatible_{\Piece_1}$.
\end{proof}

For $\relsys{X}\in \Lifts$ we denote by $W(\relsys{X})$ one of the structures
$\relsys{A}\in\Forb(\F)$ such that the structure $\relsys{X}$ is induced on $X$
by $L(\relsys{A})$. $W(\relsys{X})$ is called a {\em witness} of the fact that
$\relsys{X}$ belongs to $\Lifts$.  Note that in this definition, the witness of
a finite lift may be infinite structure, because $\Forb(\F)$ contains infinite
structures.
\begin{proof}[Proof of Theorem~\ref{mainthm}]
Clearly it suffices to prove the second part of the theorem.  By definition the
class $\Age(\Lifts)$ is hereditary, iso\-morphism\--closed, and has the joint
embedding property.  Assuming that  $\Age(\Lifts)$ has the amalgamation
property (with restrictions described), the rest of the theorem follows from the
\Fraisse{} Theorem and the fact that $\Lifts$ is the class of all lifts younger
than the \Fraisse{} limit $\relsys{U}'$ of $\Age(\Lifts)$ and thus
$\relsys{U}'$ is generic for $\Lifts$.

We show the amalgamation property.
Consider $\relsys{X},\relsys{Y},\relsys{Z}\in \Age(\Lifts)$. Assume that structure
$\relsys{Z}$ is a substructure induced by both $\relsys{X}$ and
$\relsys{Y}$ on $Z$ and without loss of generality assume that $X\cap Y=Z$.

Put
$$\relsys{A}=W(\relsys{X}),$$
$$\relsys{B}=W(\relsys{Y}),$$
$$\relsys{C}=\sh(\relsys{Z}).$$

Because $\Age(\Lifts)$ is closed under isomorphism, we can assume that
$\relsys{A}$ and $\relsys{B}$ are vertex-disjoint with the exception of
vertices of $\relsys{C}$.

Let $\relsys{D}$ be the free amalgamation of $\relsys{A}$ and $\relsys{B}$ over
vertices of $\relsys{C}$: the vertices of $\relsys{D}$ are $A\cup B$ and there
is $\vec{v}\in \rel{D}{i}$ if and only if $\vec{v}\in \rel{A}{i}$ or
$\vec{v}\in \rel{B}{i}$.

\begin{figure}
\centerline{\includegraphics{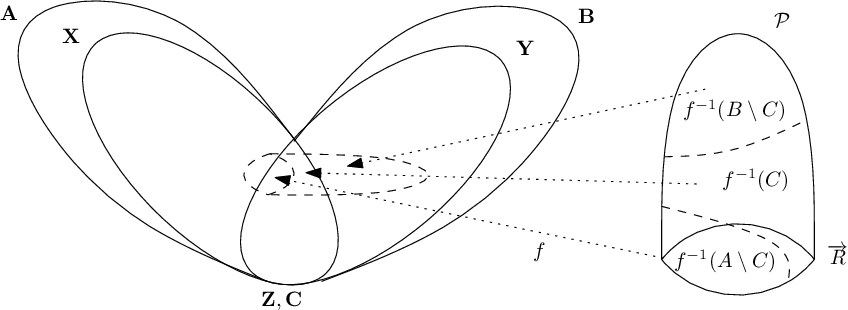}}
\caption{Construction of an amalgamation.}
\label{constructionh}
\end{figure}
We claim that the structure $$\relsys{V}=L(\relsys{D})$$ is a strong
amalgamation of $L(\relsys{A})$ and $L(\relsys{B})$ over $\relsys{Z}$ and thus
also an amalgamation of $\relsys{X}, \relsys{Y}$ over $\relsys{Z}$. The
situation is depicted in Figure~\ref{constructionh}.

First we show that the substructure induced by $\relsys{V}$ on $A$ is
$L(\relsys{A})$ and that the substructure induced by $\relsys{V}$ on $B$ is
$L(\relsys{B})$. In the other words, no new tuples to $L(\relsys{A})$ or
$L(\relsys{B})$ (and thus none to $\relsys{X}$ or $\relsys{Y}$ either) have
been introduced.
Assume to the contrary that there is a new tuple $(v_1,\ldots,v_t)\in
\ext{V}{k}$.
 By symmetry we can assume that $\{v_1,\ldots,v_t\}\subseteq A$. Explicitly, we assume that there is a piece $\Piece=(\relsys{P},\overrightarrow{R})\in E_k$ and a homomorphism $f$
from $\relsys{P}$ to $\relsys{D}$ such that
$f(\overrightarrow{R})=(v_1,v_2,\ldots, v_t)\notin
\extl{L(\relsys{A})}{k}$.

The set of vertices of $\relsys{P}$ mapped to $L(\relsys{A})$, $f^{-1}(A)$, is
nonempty, because it contains all vertices of $\overrightarrow{R}$. The set
$f^{-1}(B\setminus C)$ is nonempty $f$ is not homomorphism from
$\relsys{P}$ to $\relsys{A}$ (otherwise we would have $(v_1,v_2,\ldots, v_t)\in
\extl{L(\relsys{A})}{k})$.  Because there are no tuples spanning both vertices
$A\setminus C$ and vertices $B\setminus C$ in $\relsys D$, and because pieces
are connected, we also have $f^{-1}(C)$ nonempty. 

\begin{figure}
\centerline{\includegraphics{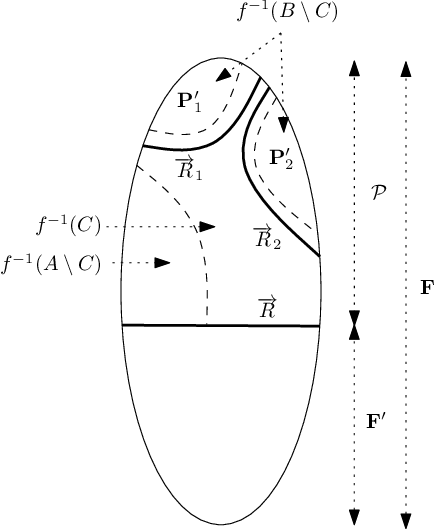}}
\caption{The decomposition of piece $\Piece$.}
\label{piecedep}
\end{figure}
We will reason about the decomposition of $\Piece$ given by $f^{-1}(C)$
and create subpieces containing vertices of $f^{-1}(B\setminus C)$.  This requires some careful analysis. The process is depicted in Figure~\ref{piecedep}.
Denote by $\relsys{F}\in \F$ the structure such that $\Piece$ is
a piece of $\relsys{F}$.
The vertices of $f^{-1}(C)$
form a g-cut in $\relsys{F}$ g-separating any vertex in $f^{-1}(B\setminus C)$ from
any vertex in $f^{-1}(A\setminus C)$ (if such a vertex exists) as well as  any vertex of $\relsys{F}\setminus P$.

We further strengthen our assumption on the choice of counter-example
(consisting of the choice of $\relsys{X}$, $\relsys{Y}$, $\relsys{A}$, $\relsys{B}$, the piece $\Piece$ and
the homomorphism~$f$):
\begin{itemize}
 \item[$(a)$] No subpiece $\Piece'=(\relsys{P}',\overrightarrow{R}')$ of $\Piece$ is a counter-example with homomorphism $f$. 
More precisely
if $f(R')\subseteq A$, then $f(\overrightarrow{R'})\in \extl{L(\relsys{A})}{i}$ for $i$ such that $\Piece'\in E_i$. Similarly
if $f(R')\subseteq B$, then $f(\overrightarrow{R'})\in \extl{L(\relsys{B})}{i}$ for $i$ such that $\Piece'\in E_i$.
 \item[$(b)$] If there is some g-component $\relsys{P}'$ of $\relsys{F}$ with g-cut $f^{-1}(C)$ contained in $f^{-1}(B\setminus C)$ such that $R\subseteq N_\relsys{P}(P')$,
then there exists a g-component $\relsys{P}''$ of $\relsys{F}$ with g-cut $f^{-1}(C)$ contained in $f^{-1}(A\setminus C)$ such that $R\subseteq N_\relsys{P}(P'')$.
\end{itemize}
It easily follows that the existence of counter-example implies the existence of counter-example
satisfying $(a)$ and $(b)$.  $(a)$ can be made to hold by considering the smallest subpiece of
$\Piece$ that is still a counter-example.  If $(b)$ fails
we can exchange $\relsys{A}$ and $\relsys{B}$ as well as exchange $\relsys{X}$ and $\relsys{Y}$.
This is possible because $f^{-1}(C)$ g-separates any vertex in $f^{-1}(A\setminus C)$ and any vertex in $f^{-1}(B\setminus C)$
and thus the existence of $\relsys{P}'$ implies $R\subseteq f^{-1}(C)$.

Denote by $\relsys{F}'$ a substructure of $\relsys{F}$ so that $R$ is a minimal
g-separating g-cut of $\relsys{F}'$ and $\relsys{P}\setminus R$ in $\relsys{F}$.
Denote by $\relsys P'_1, \relsys P'_2,\ldots, \relsys P'_l$ the substructures induced on $\relsys{P}$ by all connected
components of  $G_\relsys P\setminus f^{-1}(A)$.  We aim to find, for every $i=1,2,\ldots l$, a minimal g-separating g-cut $R_i\subseteq f^{-1}(C)$ of $\relsys{F}$ that
separates $\relsys{P}'_i$ from $\relsys{F}'$. Moreover, we will choose $R_i$ so that $R_i\not\subseteq R$. This implies the existence
of $\Piece_i=(\relsys{P}_i,\overrightarrow{R}_i)$ that is a subpiece of $\Piece$
containing $\relsys{P}'_i$.

We consider two cases:
\begin{enumerate}
\item $R\not\subseteq N_\relsys{P}(P'_i)$: Construct $R_i\subseteq N_\relsys{P}(P'_i)$ as
a minimal g-separating g-cut that g-separates $\relsys{F}'$ and $\relsys{P}'_i$ in $\relsys{F}$
(given by Proposition~\ref{prop:sep} for structure $\relsys{F}$ and g-cut $f^{-1}(C)$). 

Since pieces are connected and $R$ is a minimal g-separating g-cut in $\relsys{F}$ for $\relsys{P}\setminus R$ and $\relsys{F}'$,
we know that $\overrightarrow{R}_i$ must contain some vertex $v\notin R$.
\item $R\subseteq N_\relsys{P}(P'_i)$: In this case consider structure $\relsys{P}''$
given by $(b)$.
 Construct $R_i$ as
a minimal g-separating g-cut that g-separates $\relsys{P}''$ and $\relsys{P}'_i$ in $\relsys{F}$.

We show that $R_i\supset R$.  Because every vertex $v\in R$ is connected in $G_\relsys{P}$ to a vertex in $P_i$ and a vertex in $P''$,
we have $R_i\supseteq R$.  Moreover because pieces with roots removed
are g-components, $R$ does not g-separate $\relsys{P}''$ and $\relsys{P}'_i$ and thus $R_i\supset R$.

Since $R_i\supset R$, $R_i$ g-separates $\relsys{P}'_i$ and $\relsys{F}'$ in $\relsys{F}$.
\end{enumerate}

We have constructed a family of subpieces $\Piece_1, \Piece_2,\ldots, \Piece_l$ such that $\relsys{P}'_i$ is contained in $\relsys{P}_i\setminus \overrightarrow{R}$.
It is possible that $\Piece_i$ is a subpiece of $\Piece_j$ for some $i\neq j$.
Without loss of generality assume that $\Piece_1,\Piece_2,\ldots \Piece_{l'}$
is the maximal subset of pieces $\Piece_1, \Piece_2,\ldots, \Piece_l$ such that no piece is a subpiece of any other.
 Obviously $P\setminus f^{-1}(A)=\cup_{i=1,2,\ldots l}P'_i$ is a subset of
$\cup_{i=1,2,\ldots l'}(P_i\setminus R_i)$.

 Let $e_i$ be the index of the equivalence class
of $\sim$ such that $\Piece_i\in E_{e_i}$.
Now we use assumption $(a)$. All the pieces $\Piece_i$,
$i=1,\ldots, l'$ are subpieces of $\Piece$.  Thus we have that $f(\overrightarrow{R}_i)\in
\extl{L(\relsys{D})}{e_i}\implies f(\overrightarrow{R}_i)\in \extl{L(\relsys{A})}{e_i}$.
Thus there exists a piece $\Piece^A_i=(\relsys{P}^A_i,
\overrightarrow{R}^A_i)$, $\Piece^A_i\sim \Piece_i$, and a homomorphism
$f^A_i$ from $\relsys{P}^A_i$ to $\relsys A$ such that
$f^A_i(\overrightarrow {R}^A_i) = f (\overrightarrow {R_i})$, for
every $i=1,2,\ldots, l'$.

In this situation we want to create $\Piece^A=(\relsys{P}^A,\overrightarrow{R}^A)$ as a
copy of $\Piece$ with $\Piece_i$ replaced by $\Piece^A_i$ for every $i=1,\ldots
l'$.  By (repeated) application of Lemma~\ref{eqkousky} we will then have $\Piece^A\sim
\Piece$. 
To make this possible, we must show that no root vertex $v$ of $\Piece_i$ is
contained in some $\relsys{P}_j\setminus R_j$ for $1\leq i\leq l'$, $1\leq
j\leq l'$.  (Otherwise replacing $\Piece_j$ by $\Piece^A_j$ may make it
impossible to replace $\Piece_i$ by $\Piece^A_i$.)  Assume, to the contrary
that there is such a choice of $\Piece_i$, $\Piece_j$ and $v$. Because $\Piece_i$
is not subpiece of $\Piece_j$, nor vice versa, there is a root $v'$ of $\Piece_i$ that
is not contained in $\relsys{P}_j$.  Because $v$ is in $\relsys{P}_j\setminus R_j$ that is a g-component of $\relsys{F}$ with cut $R_j$  and $v'$ is not, we conclude that $v$ and $v'$  are g-separated by $R_j$ in $\relsys{F}$. This leads to the fact that
$R_j\cap P_i$ is g-cut in $\relsys{P}_i$ g-separating $v$ and $v'$. This is not possible
because, by our construction, $v,v'\in N_\relsys{P}(P'_i)$ and $\relsys{P}'_i$ is a g-component of $\relsys{P}$ with g-cut $f^{-1}(C)$. Thus $v$ and $v'$ are connected by a walk in $G_{\relsys{P}_i}$ containing only vertices of $\relsys{P}'_i$. It is not possible for vertex of $\relsys{P}'_i$ to be in $R_j$, because vertices of $R_j$ are in $f^{-1}(C)$, while vertices of $\relsys{P}'_i$ are in $f^{-1}(B\setminus C)$.

Finally define $f^A:P^A\to A$ as follows:
\begin{enumerate}
 \item $f^A(x)=f^A_i(x)$ when $x\in P^A_i$ for some $i=1,2,\ldots, l'$;
 \item $f^A(x)=f(x)$ otherwise.
\end{enumerate}
It is easy to see that $f^A$ is a homomorphism from $\relsys P^A$ to
$\relsys{A}$ such that $f^A(\overrightarrow{R}^A)=(v_1,v_2,\ldots
v_t)$. This is a contradiction with $(v_1,v_2,\ldots v_t)\notin
\extl{L(\relsys{A})}{k}$.

It remains to verify that $\relsys{D}\in \Forb(\F)$. We proceed analogously.
Assume that $f$ is a homomorphism from some $\relsys{F}\in \F$ to $\relsys{D}$
and further assume that the counter-example is chosen in a way so
$\relsys{F}\setminus f^{-1}(A)$ has minimal number of vertices.

Because $\relsys{A},\relsys{B}\in \Forb(\F)$, $f$ must use vertices of
$\relsys{D}\setminus A$ and vertices of $\relsys{D}\setminus B$ and, because $\relsys{F}$ is connected, also vertices of $\relsys{C}$.
Analogously to the previous part, $f^{-1}(C)$ forms a g-cut in $\relsys{F}$. Denote by
$\relsys{F}_A$ a g-component of $\relsys{F}$ with g-cut $f^{-1}(C)$ contained in $\relsys{F}\setminus f^{-1}(B)$ and by $\relsys{F}_B$ a
g-component of $\relsys{F}$ with g-cut $f^{-1}(C)$ contained
in $\relsys{F}\setminus f^{-1}(A)$.
Denote by $R$ a minimal g-separating g-cut in
$\relsys{F}$ contained in $f^{-1}(C)$ that g-separate $\relsys{F}_A$ and $\relsys{F}_B$ (given by Proposition~\ref{prop:sep}). 
Denote by $\Piece=(\relsys{P},\overrightarrow{R})$ a piece of $\relsys{F}$ containing $\relsys{F}_B$.
We have shown that $f(\overrightarrow{R})\in \extl{L(\relsys{A})}{i}$, for $i$ such that $\Piece\in E_i$.
Consequently there is $\Piece'=(\relsys{P}',\overrightarrow{R}')$, $\Piece\sim \Piece'$, and a homomorphism $f^A:\relsys{P}'\to \relsys{A}$
such that $f^A(\overrightarrow{R}')=f(\overrightarrow{R})$.
Now denote by $\relsys{F}'$ a structure created from $\relsys{F}$ by replacing piece $\Piece$
by piece $\Piece'$ and identifying $\overrightarrow{R}$ with $\overrightarrow{R}'$. Because $\Piece\sim\Piece'$, $\relsys{F}'$ is isomorphic to some structure in $\F$. Now consider the homomorphism $f':\relsys{F}'\to \relsys{D}$ defined as follows:
\begin{enumerate}
\item $f'(x)=f^A(x)$ for $x\in P'$,
\item $f'(x)=f(x)$ otherwise.
\end{enumerate}
The size of $\relsys{F}'\setminus f'^{-1}(A)$ is strictly smaller than the size of $\relsys{F}\setminus f^{-1}(A)$, a contradiction with the minimality of the counter-example.

This finishes the proof of the amalgamation property of $\Age(\Lifts)$: while
$\relsys{V}$ may be infinite (because witnesses may be infinite), the lift
$\relsys{V}'$ induced on vertices $X\cup Y$ by $\relsys{V}$ is 
the finite amalgamation of $\relsys{X}, \relsys{Y}$ over $\relsys{Z}$ and thus
$\relsys{V}'\in \Age(\Lifts)$.
\end{proof}

\section{Non-existence of universal structures for classes $\Forb(\F)$.}
\label{sec:classification}

In this section we show that there exists infinite families $\F$ such that
there is no countable $\omega$-categorical universal structure for $\Forb(\F)$.
This is in contrast with the finite case, where the universal structures always
exists.
\begin{thm}
\label{mainthm2}
Let $\F$ be a family of finite connected relational structures (of finite type).
Assume that:
\begin{enumerate}
 \item[$(i)$] The size of all minimal g-separating g-cuts of structures in $\F$ is bounded by $n$.
 \item[$(ii)$] 
Let $\Piece=(\relsys{P},\overrightarrow{R})$ and $\Piece'=(\relsys{P}',\overrightarrow{R}')$ be two pieces (of some structures in $\F$).
Denote by $\APiece$ a rooted structure such that $\Piece\oplus \APiece=\relsys{F}\in \F$.  If $\Piece'\oplus \APiece$ is defined and $\Piece'\oplus \APiece\notin \Forb(\F)$, then there is $\relsys{F}'\in \F$ isomorphic to $\Piece'\oplus \APiece$.
\end{enumerate}
Then the following conditions are equivalent:
\begin{enumerate}
 \item[$(a)$] $\F$ is a regular family of connected structures.
 \item[$(b)$] There is an ultrahomogeneous lift $\extsys{U}$ with only finitely many new relations such that its shadow $\sh(\extsys{U})$ is a universal structure for the class $\Forb(\F)$.
 \item[$(c)$] There exists an countable $\omega$-categorical universal structure for $\Forb(\F)$.
\end{enumerate}
\end{thm}
\begin{proof}
$(a)\implies (b)$ follows from Theorem~\ref{mainthm} for the class $\F$.

$(b)\implies (c)$ is immediate. The shadow of every ultrahomogeneous structure
with finitely many relations is $\omega$-categorical.

To see that $(c)\implies (a)$, assume to the contrary the existence of
$\F$ satisfying $(i)$ and $(ii)$ which is not regular such that there is a universal structure
$\relsys U\in \Forb(\F)$ which is $\omega$-categorical. 

Because the sizes of minimal g-separating g-cuts are bounded by $n$, we know that
there is $n'\leq n$ with infinitely many pieces $\Piece_1,\Piece_2,\ldots$ of
width $n'$ such that the corresponding sets
$\Incompatible_{\Piece_1},\Incompatible_{\Piece_2},\ldots$ are all different.

 From Theorem~\ref{Nardzewsky} it follows that there are only
finitely many orbits of $n'$-tuples.  Denote by $k$ the number of orbits of $n'$ tuples.
Now assign every piece $\Piece_i$ the set $O_i$ of all orbits such
that there exists a rooted embedding from $\Piece_i$ to $\relsys{U}$ (that is a rooted homomorphism that is also an embedding of $\relsys{P}_i$ to $\relsys{U}$) sending
the root of $\Piece_i$ to the orbit.  All the sets $O_i$ are finite of size at most $k$. By the
pigeonhole principle there is $i\neq j$ such that $O_i=O_j$.  

By our assumption $(ii)$ we know that for two pieces $\Piece=(\relsys{P},\overrightarrow{R})$ and $\Piece'=(\relsys{P}',\overrightarrow{R}')$ such that $\relsys{P}, \relsys{P}'\notin \Forb(\F)$ we have $\Incompatible_\Piece = \Incompatible_{\Piece'}$ if and only if the rooted structure induced by $\relsys{P}$ on $\overrightarrow{R}$ is identical
to the rooted structure induced by $\relsys{P}'$ on $\overrightarrow{R}'$.
Consequently all those pieces belong to the same class of $\sim$ and there are only finitely many classes $\sim$ containing piece $\Piece=(\relsys{P},\overrightarrow{R})$ such that $\relsys{P}\notin \Forb(\F)$.
We can thus assume that $\relsys{P}_i, \relsys{P}_j\in \Forb(\F)$ and thus there is an isomorphic copy of both $\relsys{P}_i$ and $\relsys{P}_j$ in
$\relsys{U}$ (we where choosing $i$ and $j$ from infinitely many equivalence classes of $\sim$). 

Now, because
$\Incompatible_{\Piece_i}\neq \Incompatible_{\Piece_j}$ there is a rooted structure $\APiece$ that distinguishes
$\Incompatible_{\Piece_i}$ from $\Incompatible_{\Piece_j}$.  Without loss of generality assume
that $\APiece\in \Incompatible_{\Piece_i}$.  By our assumption $(ii)$, $\APiece\oplus \Piece_j\notin
\F$ implies $\APiece\oplus \Piece_j \in \Forb(\F)$.  Consequently
there is an embedding from $\APiece\oplus \Piece_j$ to $\relsys{U}$. This
embedding must map the root of $\APiece\oplus \Piece_j$ to a tuple within an
orbit $o\in O_j$.  Since $\APiece\oplus \Piece_i\notin \Forb(\F)$, we also have
$o\notin O_i$. This is in contradiction to $O_i=O_j$.
\end{proof}

\paragraph{Example.}
Consider relational structures with two binary relations (red and blue directed
edges).  The family consist $\F_\mathrm{rb-balanced}$ of all red-blue balanced
cycles, defined as follows.  Take two paths $P_1$ and $P_2$ with initial and terminal
vertex specified where $P_1$ consists only of red edges, $P_2$ consists
only of blue edges, and $P_1$, $P_2$ have have identical non-zero
algebraic length. (The algebraic length of a path is the number of forwarding
edges minus number of backwarding edges on the way from the initial vertex to
the terminal vertex.) A {\em red-blue balanced cycle} is created as a disjoint
union of $P_1$ and $P_2$ with the corresponding initial and terminal vertices
identified.  

The pieces of $\F_\mathrm{rb-balanced}$ consist of paths starting with a (possibly
empty) blue segment followed by a (possibly empty) red segment.
The equivalence class of a piece can be characterized by
\begin{enumerate}
\item the existence of red segment and the existence of blue segment,
\item whether the first root is in red or blue segment,
\item the algebraic length of red segment (if any) and the algebraic length of blue segment (if any).
\end{enumerate} 
It is not difficult to see that decision whether $\Piece\oplus\APiece$ is in
$\Forb(\F_\mathrm{rb-balanced})$ depends precisely on the information above.  There
are infinitely many different algebraic lengths and thus also equivalence
classes of $\sim$. Consequently $\F_\mathrm{rb-balanced}$ is not a regular
family.

Finally we verify that $\F_\mathrm{rb-balanced}$ satisfies the assumptions of
Theorem~\ref{mainthm2}. Cycles have minimal separating cuts of size 2 giving
$(i)$. We now verify $(ii)$. The only rooted structures $\APiece$ such that
$\Piece\oplus\APiece$ is isomorphic to some red-blue balanced cycle are again
pieces of red-blue balanced cycles. The algebraic lengths are preserved in
homomorphic images and because the algebraic lengths are non-zero, no red-blue
balanced cycle can be mapped to a path. Consequently for every piece of
red-blue balanced cycle $\Piece$, $\Piece'\oplus\APiece$ is a cycle and thus
$\Piece'\oplus\APiece\notin \Forb(\F_\mathrm{rb-balanced})$ imply that
$\Piece'\oplus\APiece$ is a red-blue balanced cycle giving $(ii)$.
Consequently, by Theorem~\ref{mainthm2}, there is no $\omega$-categorical universal structure for
$\Forb(\F_\mathrm{rb-balanced})$.

There are non-regular families $\F$ where universal graph for $\Forb(\F)$
exists.  Consider relational structures with one binary relation. For
simplicity, assume that loops and bi-directional edges are forbidden so we
consider oriented graphs.  Family $\F_\mathrm{balanced}$ of all balanced
orientations of graph cycles (i.e., orientations having the same number of
forward and backward edges) is not a regular family.  Here all pieces are all
oriented paths.  The equivalence class of a piece depends on the algebraic
length of the path.
On the other hand, however, homomorphic image of balanced 4-cycle is an
oriented edge. Consequently the countably infinite graph with no edges is the
universal graph for $\Forb(\F_\mathrm{balanced})$.  

The family $\F_\mathrm{balanced}$ do not satisfy requirement of
$Theorem~\ref{mainthm2}$, $(ii)$. Let $\Piece$ be a piece of balanced cycle and
$\APiece$ an oriented path rooted in the endpoints.  $\Piece\oplus\APiece$
is a oriented cycle and because it contains an edge we know that
$\Piece\oplus\APiece \notin \Forb(\F)$. Every oriented cycle can be
constructed this way. The smallest family $\F'$ containing
$\F_\mathrm{balanced}$ that satisfy $(ii)$ is thus the family of all oriented cycles.
Such $\F'$ is a regular family.

Further applications of Theorem~\ref{mainthm2} are given in the following section.

\section{Homomorphism dualities and constraint satisfaction problems}
\label{sec:dualities}

A constraint satisfaction problem (CSP) is the following decision problem:

\smallskip

\noindent
{\bf Instance:} A finite structure $\relsys{A}$.

\noindent
{\bf Question:} Does there exist a homomorphism $\relsys{A}\to\relsys{D}$?
\smallskip

\noindent We denote by $\CSP({\mathcal D})$ the class of all finite structures
$\relsys{A}$ with $\relsys{A}\to\relsys{D}$ for some $\relsys{D}\in {\mathcal
D}$.

Recall that a {\em homomorphism duality} (for structures of given type) is any
equation $$\Forb(\F)=\CSP({\mathcal D}).$$ When both $\F$ and $\mathcal D$ are
finite sets of finite structures, we call the pair $(\F, \mathcal D)$ {\em a
finite duality pair} \cite{NesetrilPultr,NTardif,Hell}. When $\F$ is an
infinite set of finite structures, and $\mathcal D$ is a finite set of finite
structures, we call it an {\em infinite-finite duality} \cite{Tardif}.

Dualities play a role not only in complexity problems but also in logic, model
theory, the theory of partial orders and category theory. In particular, it
follows from \cite{Atserias} and \cite{Rossman} that dualities coincide with
those first-order definable classes which are homomorphism-closed. 

For the sake of simplicity, in the following discussion we shall restrict
ourselves to the case where $\mathcal{D}$ consists of a single element
$\relsys{D}$.  $\relsys{D}$ is called the {\em dual of $\F$} (it is easy to see
that $\relsys{D}$ is up to homomorphism-equivalence uniquely determined). 

The notion of universal structures and duals is related.  Given a class $\K$ of
countable structures, an object $\relsys{U}\in\K$ is called {\em hom-universal}
for $\K$ if for every object $\relsys{A}\in \K$ there exists a homomorphism
$\relsys{A}\to \relsys{U}$. The following is immediate from the definitions:
\begin{prop}
Let $\F$ be a family of relational structures. Structure $\relsys{D}$ is the
dual of $\F$ if and only if $\relsys{D}$ is hom-universal for $\Forb(\F)$. 
\end{prop}
In this section we shall show how to turn the universal structure constructed
in Section~\ref{sec:mainthm} into a finite dual.
This is possible only in the special cases where a finite dual exists.  First
we review some results characterizing dualities.

A (relational) tree can be defined as follows (see \cite{NTardif}): The {\em incidence graph}
$\ig(\relsys{A})$ of relational structure $\relsys{A}$ is the bipartite graph
with parts $A$ and $\Block(A)$, where
$$\Block(A) = \{ (i,(a_1, \ldots, a_{\delta_i})) : i \in I, (a_1, \ldots,
a_{\delta_i}) \in \rel{A}{i} \},$$
and edges $[a, (i,(a_1, \ldots, a_{\delta_i}))]$ such that                      
$a \in (a_1, \ldots, a_{\delta_i})$.
(Here we write $x \in (x_1, \ldots, x_n)$ when there exists an index $k$
such that $x = x_k$; $\Block(A)$ is a multigraph.) Relational structure
$\relsys{A}$ is called a {\em (relational) tree} when
$\ig(\relsys{A})$ is a graph tree (see e.g. \cite{MatN}).  The definition of
relational trees by the incidence graph $\ig(\relsys{A})$ allows us to use
graph terminology for relational trees. 

\begin{thm}[\cite{NTardif}]
\label{dualpairs}
For every finite family $\F$ of finite relational trees there exists a dual
$\relsys{D}$. Up to homomorphism-equivalence there are no other
finite dualities with only one dual.
\end{thm}
Various constructions of duals of a given $\F$ are known
\cite{NTardif2}. 
 More recently, infinite-finite dualities have been
characterized:
\begin{thm}[\cite{Tardif}]

\label{thm:tardif1}
All regular families $\F$ of relational trees have a finite dual $\relsys{D}$.
\end{thm}

\begin{thm}[\cite{Tardif}]
\label{thm:tardif2}
The family $\F$ of relational trees has a finite dual if and only if its upward
closure $\UP(\F)$ is regular.
\end{thm}

Here the {\em upward closure}, $\UP(\F)$, is the class of all relational trees
$\relsys{T}_1$ such that there is $\relsys{T}_2\in \F$ and $\relsys{T}_2\to
\relsys{T}_1$.

We remark that all these characterizations extend naturally to duality pairs
$(\F,\mathcal D)$ where structures in the class $\F$ are not necessarily
connected (i.e., they are relational forests). In this case however $\mathcal D$
generally consists of one or more structures.  See \cite{Foniok, Tardif} for
details.

The construction of Section~\ref{sec:mainthm} may be used to obtain an
alternative way of constructing a dual in the proof of Theorems~\ref{dualpairs}
and~\ref{thm:tardif1}:

\begin{corollary}[of Theorem~\ref{mainthm}]
\label{nasdual}
Let $\F$ be a regular set of finite relational trees. Then there exists a class
$\Lifts$ of monadic lifts such that:
\begin{itemize}
  \item[$(i)$] $\Age(\Lifts)$ is an amalgamation class with free amalgamation;
  \item[$(ii)$] The \Fraisse{} limit of $\Age(\Lifts)$ is an ultrahomogeneous structure $\relsys{U}'$ such that $\sh(\relsys{U}')=\relsys{U}$ is universal for $\Forb(\F)$;
  \item[$(iii)$] $\relsys{U}'$ has a finite retract $\extsys{D}$ and consequently $\sh(\extsys{D})=\relsys{D}$ is a dual of $\F$.
\end{itemize}
\end{corollary}

\begin{proof}[Proof]
Observe that the minimal g-separating g-cuts of a relational tree all have size 1.
Thus for a fixed family $\F$ of finite relational trees Theorem~\ref{mainthm}
establishes the existence of a class $\Lifts$ and lift $\relsys{U}'$ satisfying
$(ii)$.  Class $\Age(\Lifts)$ is closed under strong amalgams that are free in
the shadow.  With only unary relations added to the lift, we immediately get
that $\Age(\Lifts)$ is closed under free amalgamation, too, thereby obtaining
$(ii)$.

We show $(iii)$.  We find finite $\relsys{D}'$ which is a retract of
$\relsys{U}'$, $\sh(\relsys{D}')\in \Forb(\F)$, and for which there is a
homomorphism $\relsys{A}\to \sh(\relsys{D}')$ if and only if there is a
homomorphism $\relsys{A} \to \sh(\relsys{U}')$ for every relational structure
$\relsys{A}$. 

Construct $\relsys{D}'$ from $\relsys{U}'$ by identifying all vertices of the
same color (recall that the  color of vertex $v$ is the set $\{i;(v)\in
\ext{U}{i}\}$).  Denote by $r$ the homomorphism (retraction) $\relsys{U}'\to
\relsys{D}'$.  Obviously, if $f:\relsys{Y}\to \relsys{U}'$ is a homomorphism
then $f\circ r:\relsys{Y}\to \relsys{D}'$ is also a homomorphism. Thus
$\relsys{A} \to \sh(\relsys{U}')$ implies $\relsys{A}\to \sh(\relsys{D}')$.

It remains to show that $\sh(\relsys{D}')\in \Forb(\F)$. Suppose, to the
contrary, that there is a tree $\relsys{F}\in \F$ and a homomorphism
$f:\relsys{F}\to \sh(\relsys{D}')$.  Let $\relsys{X}$ be a lift created from
$\relsys{F}$ by adding an extended relation $(v)\in \ext{X}{i}$ if and only if
$(f(v))\in \ext{D'}{i}$, for every $i=1,2,\ldots N$. Obviously $f$ is also a homomorphism
$\relsys{X}\to \relsys{D}'$.
Consider lift $\relsys{Y}$ induced by $\relsys{X}$ on elements of some tuple
$\vec{v}\in \rel{X}{j}$.  Since $\relsys{F}$ is a relational tree, the shadow
$\sh(\relsys{Y})$ has only one tuple. Because $\relsys{D}'$ is retract of
$\relsys{U}'$ we know that the homomorphic image of $\relsys{Y}$ is a substructure 
of $\relsys{U}'$ and thus it is in $\Lifts$.

Because $\relsys{F}$ is a relational tree, it is possible to construct a homomorphic
copy of $\relsys{X}$ by starting with the homomorphic image of $\relsys{Y}$ in
$\Lifts$ and using free amalgamation (over a one-element set) to add lifts of
homomorphic images of all other tuples of $\relsys{F}$.  It follows that the
homomorphic image of $\relsys{X}$ is in $\Lifts$, a contradiction.
\end{proof}

\paragraph{Remark.}
We stress the fact that families of trees are not the only regular families
$\F$ of relational structures where the universal structure for $\Forb(\F)$ can
be described as a shadow of an ultrahomogeneous monadic lift $\relsys{U}'$.
For example, consider relational structures created from a relational tree by
replacing tuples by an arbitrary irreducible structure (recall that a structure
is {\em irreducible} if it has no vertex cuts). Such structures have all
minimal g-separating g-cuts of size 1. One can easily construct continuum many
such examples. There is however no finite retract of $\relsys{U}'$ satisfying
the statement of Corollary~\ref{nasdual}.  (Such structures
cannot be constructed from individual tuples by the aid of free amalgamation.)
Of course such structures may have infinite chromatic numbers.

Note that it is also possible to construct the dual $\relsys{D}$ of $\F$
without using the \Fraisse{} limit. Also this follows by our construction in
Theorem~\ref{mainthm}. For every possible combination of new relations on a
single vertex, create a single vertex of $\relsys{D}$ and then keep adding
tuples as long as possible so that $\relsys{D}$ is still in $\Lifts$ (in a
similar way to the proof of Proposition~\ref{prop:duals}).

We have shown that special cases of universal structures can be used to
construct duals.  Now we show the opposite: every dual can be turned into a
universal structure by an especially simple monadic lift.

\begin{prop}
\label{prop:duals}
For a family $\F$ of relational structures the following two statements are equivalent:
\begin{enumerate}
 \item[$(i)$] There is a finite dual $\relsys{D}$ of $\F$.
 \item[$(ii)$] There exists a finite family $\F'$ of monadic lifts $\relsys{X}$ whose shadow $\sh(\relsys{X})$ has one tuple with the following property:

Denote by $\Lifts$ the class of all lifts $\relsys{Y}$ such that
    \begin{enumerate}
     \item[(a)] $\relsys{Y}\in \Forb(F')$, and
     \item[(b)] every vertex of $\relsys{Y}$ is in precisely one extended relation $\ext{Y}{i}$.
    \end{enumerate}
   There is a generic lift $\relsys{U}'$ for $\Lifts$ and its shadow $\sh(\relsys{U}')$ is an $\omega$-categorical universal structure for $\Forb(\F)$.
\end{enumerate}
\end{prop}
Loosely speaking, the class $\Lifts$ is described by forbidden colors of vertices and forbidden colorings of edges.

\begin{proof}
$(i)\implies (ii)$:
Fix a dual $\relsys{D}$ with vertices $\{1,2,\ldots N\}$ and consider lifts
with $N$ extended unary relations.  Let $\F'$ be the family of all structures
$\relsys{X}$ such that:
\begin{enumerate}
\item the vertex set of $\relsys{X}$ is $X\subseteq \{1,2,\ldots N\}$;
\item there is a tuple $\vec{v}\in \rel{X}{j}$ for some $j\in I$ such that $\vec{v}\notin \rel{D}{j}$;
\item for every $i\in X$ there is a tuple $(i)\in \ext{X}{i}$;
\item there are no other tuples in $\relsys{X}$ and there are no vertices in $X$ except ones in $\vec{v}$.
\end{enumerate}
By definition, $\Age(\Lifts)$ is obviously an (free) amalgamation class (all forbidden
substructures are irreducible). 

We show that the shadow of $\Lifts$ is $\Forb(\F)$.  For every $\relsys{A}\in
\Forb(\F)$ and homomorphism $f:\relsys{A}\to \relsys{D}$ construct a lift
$\relsys{X}$ by putting $(v)\in \ext{X}{i}$ if and only if $f(v)=i$. It is easy
to see that $\relsys{X}\in \Lifts$. On the other hand, for every structure
$\relsys{A}$ and lift $\relsys{X}\in \Lifts$, a homomorphism $\relsys{A}\to
\sh(\relsys{X})$ can be interpreted as an $\relsys{D}$-coloring of $\relsys{A}$
and thus $\relsys{A}\in \Forb(\F)$.

The rest of statement follows by \Fraisse{} theorem analogously as Theorem~\ref{mainthm}.

In the opposite direction assume the existence of $\F'$, $\Lifts$ and $\relsys{U}'$ satisfying the
statement of the proposition. Construct the retract
$\relsys{D}'$ of $\relsys{U}'$ by unifying all
vertices of the same color.   This gives a homomorphism (retraction) $r:\sh(\relsys{U}')\to
\relsys{D}'$. Put $\relsys{D}=\sh(\relsys{D}')$. We show that $\relsys{D}$ is the dual of $\F$.

For every $\relsys{A}\in \Forb(\F)$ there is an embedding
$e:\relsys{A}\to\sh(\relsys{U}')$.  It follows that $\relsys{A}\in
\CSP(\{\relsys{D}\})$ because $e\circ r$ is a homomorphism $\relsys{A}\to
\relsys{D}$.  To see that $\relsys{D}\notin \Forb(\F)$, assume for a
contradiction that there is $\relsys{F}\in \F$ and a homomorphism
$f:\relsys{F}\to \relsys{D}$. Create lift $\relsys{X}$ from $\relsys{F}$ by
adding a tuple $(v)$ to $\ext{X}{i}$ if and only if $(f(v))\in \ext{D'}{i}$.
Lift $\relsys{X}$ satisfy condition $(a)$ of the definition of $\Lifts$.  For
every $\relsys{F}'\in \F$ a homomorphism $\relsys{F}'\to \relsys{X}$ implies a
homomorphism $\relsys{F}'\to\relsys{U}$ giving (b) and thus $\relsys{X}\in
\Lifts$. A contradiction with $\relsys{U}'$ being generic for $\Lifts$ and
$\sh(\relsys{U}')\in \Forb(\F)$.
\end{proof}

Theorem~\ref{mainthm2} and Proposition~\ref{prop:duals} imply Theorem~\ref{thm:tardif2}. It is easy to see that $\UP(\F)$ have the size of separating
cuts o bounded by 1 and moreover it is closed for free amalgams of trees.

To show both implications of Theorem~\ref{dualpairs} it is necessary to show
the non-existence of monadic lifts as described in Proposition~\ref{prop:duals}
for families $\F$ not consisting of relational trees.  This is possible with a
more systematic study of the minimal arities needed in the lift for a given
family $\F$, and by giving a more explicit description of the lifts via
forbidden substructures, as shown in \cite{paper2}.

\section{Acknowledgment}

We would like to thank to Manuel Bodirsky for informing us about \cite{Covington}
and to Jan Foniok for pointing out the connection to the
classification of infinite-finite dualities in \cite{Tardif}.

We are also grateful to Andrew Goodall for remarks and comments that led to
improvements in the quality of this paper.

Finally we would like to thank to the anonymous referee for correcting two
principal examples in the paper and number of useful remarks.

\end{document}